\documentclass[11pt,a4paper]{article}

\usepackage[latin1]{inputenc}

\usepackage[normalem]{ulem}

\usepackage{amssymb}

\usepackage{graphicx}

\usepackage{graphicx,indentfirst,amsmath,amsfonts,amssymb,amsthm,newlfont}
\usepackage{epsfig}
\newtheorem{theorem}{Theorem}[section]

\newcommand{\Prob}{\mathbf P}
\newcommand{\R}{\mathbf R}
\newcommand{\N}{\mathbf N}

\begin{document}
\begin{center}


 {\Large {\bf Stochastic delay differential equations with jumps \\[2mm]
in differentiable manifolds
 }}

\end{center}

\vspace{0.3cm}

\begin{center}
{\large { Leandro Morgado}\footnote{E-mail:
morgado@ime.unicamp.br.
Research supported by  FAPESP 11/14797-2.} \ \ \ \ \ \ \ \ \ \ \ \ \
{ Paulo R. Ruffino}\footnote{Corresponding author, e-mail:
ruffino@ime.unicamp.br.
Partially supported by FAPESP nr. 12/18780-0, 
11/50151-0 and CNPq nr. 477861/2013-0.}}

\vspace{0.2cm}

\textit{Departamento de Matem\'{a}tica, Universidade Estadual de Campinas, \\
13.083-859- Campinas - SP, Brazil.}

\end{center}

\begin{abstract} In this article we propose a model for stochastic delay
differential equation with jumps (SDDEJ) in a differentiable manifold $M$
endowed with a connection $\nabla$. In our model, the continuous part is driven
by vector fields with a fixed delay and the jumps are assumed to come from a
distinct source of (c\`adl\`ag) noise, without delay. The jumps occur along
adopted differentiable curves with some dynamical relevance (with fictitious
time) which allow to take parallel
transport along them. Using a geometrical approach, in the last section, we
show that the horizontal lift of the solution of an SDDEJ is again a solution
of an SDDEJ in the linear frame bundle $BM$
with respect to a connection $\nabla^H$ in $BM$.

\end{abstract}


\vspace{0.3cm}
\noindent {\bf MSC2010 subject classification:} 60H10, 34K50, 53C05.


\section{Introduction}

Many natural phenomena present delays with respect to inputs: mathematically,
this is a well known and established theory in the literature. In fact,
standard delay differential equations have been
extensively studied, see e.g. the classical Hale \cite{hale} and references
therein. Models
for
these equations in manifolds appear in Oliva \cite{oliva},
and stochastic perturbations are considered in Langevin, Oliva and
Oliveira
\cite{langevin}, Mohammed \cite{mohammed}, Mohammed and Scheutzow \cite{Scheutzow, Scheutzow2},
Caraballo, Kloeden and Real \cite{Caraballo-Kloeden} for SPDE, among
many others. More recently,
differential equations with random unbounded delay have been considered in
Garrido-Atienza, Ogrowsky and Schmalfuss
\cite{Ogro-Gar-Schmal}. Besides delay, another
usual characteristic in systems in biology, physics, economics, climatology,
etc
is the presence of jumps, both in the input and in the output.
In this paper, we put these two characteristics, delays and 
jumps,
together
in the same mathematical framework.

In our model of  stochastic delay
differential equation with jumps (SDDEJ), the continuous part of
the solution is driven by vector fields with a fixed delay $d > 0$, and the
jumps are assumed to come from a
distinct source of (c\`adl\`ag) noise, without delay.
This idea is inspired by the fact that
in several phenomena,  informations reach a receptor by different
communication sources (or channels), hence it is reasonable that delays on time
are dependent on these distinct sources. As a simple example: in a storm,
lightnings
have instantaneous impact, but thunders come with delay.

Following the ideas behind the so called Marcus equation, as in Kurtz,
Pardoux and Protter \cite{KPP}, our jumps in the solution occur along
fictitious differentiable curves which allow to take parallel
transport along them. Generically speaking, as in Section 2, these fictitious
differentiable curves can be taken in many distinct ways (say, randomly or
along geodesics, etc). In our model, presented in Section 3, they follow the
deterministic flow generated by the vector fields, without delay.
We consider that the noise splits into a continuous component which is
 a
Brownian motion with drift (extensible to a class of continuous
semimartingales) plus a component given by a sequence of
(c\`adl\`ag) jumps. The number of jumps is assumed to be finite in a compact
time interval. Hence it includes L\'evy-jump diffusion, but not L\'evy process
in general. This idea of finite jumps in bounded intervals has parallel in the
theory of chain control sets, see Colonius and Kliemann 
\cite{Colonius-Kliemann}, Patr\~ao and San Martin \cite{sanmartin},
and references therein.

Another model with some numerical results for delay stochastic systems with 
jumps in Euclidean spaces can
be found in Dareiotis, Kumar and Sabanis \cite{Dar-Kumar-Sabanis} 
for L\'evy processes. Also, stochastic geometry
with jumps is considered in Cohen \cite{Cohen}, where the authors use second
order calculus.


%
%

In our approach, the delay are treated using parallel transport along the
solutions, prescribed by a connection $\nabla$ in a differentiable manifold
$M$. In Catuogno and Ruffino \cite{PP}, the authors consider a
geometrical approach to stochastic delay differential equations (SDDE) on a
manifold $M$. In particular, they prove that the horizontal lift of a SDDE
solution to the linear frame bundle $BM$ is, again, a solution of an associated
SDDE in $BM$, with respect to a horizontal connection $\nabla^H$ in $BM$.

The paper is organized as follows: In Section 2 we construct parallel
transports along curves with jumps: here, the 
jumps are
taken along a generic family of fictitious curves. In Section 3, our model of
delay differential equations with jumps is presented. Compactness of
the manifold $M$ is assumed only to guarantee the existence and uniqueness of
solution for SDDE (without jumps) for all $t\geq 0$, as in L\'eandre and
Mohammed \cite{leandre}. Finally, in Section 4, we explore geometrical aspects
of the SDDEJ. After a short revision on the geometry of linear frame bundle and
 on the horizontal connection $\nabla^H$ in $BM$, we show that
 the horizontal lift of solutions of SDDEJ in $M$ are again solutions
of SDDEJ in $BM$, extending the result mentioned above in \cite{PP}.

\section{General aspects of parallel transport}

Let $M$ be a differentiable manifold, and $\nabla$ a connection on $M$. This
structure, via parallel transport, allows one
to map vectors from a tangent space at a point in a differentiable curve into
the tangent space at another point of this curve. To fix notation, consider a
differentiable curve $\alpha: I \rightarrow M$ defined in an interval
$I\subset \R$. For
$s,t \in I$, the parallel transport along
$\alpha$ from $\alpha(s)$ to $\alpha(t)$ induced by $\nabla$
is the linear isometry denoted by $P^\nabla_{s,t}
(\alpha): T_{\alpha(s)} M \to T_{\alpha(t)} M$, such that the covariant
derivative of $t\mapsto P^\nabla_{s,t}
(\alpha)$ vanishes. If $\alpha$ is continuous and differentiable by parts, its
parallel transport is constructed joining the corresponding
parallel transports along each differentiable segment (see e.g. Kobayashi and
Numizu \cite{nomizu}).

\bigskip

\subsection{Parallel transport along a curve with jumps}

Let $\gamma: I \to M$ be a càdlàg curve with discontinuities in a
countable, discrete and closed set $D = \{t_1, t_2, \ldots \}$, possibly 
finite. Suppose that
$\gamma$ is differentiable in $ I \setminus D$. Let $\mathcal{B} =
(\beta_n)_{n \in \N}$ be a family of differentiable curves $\beta_n: [0,1] \to
M$ such that, for each $n \in \N$, $\beta_n (0) = \displaystyle \lim_{s \to
t_n-} \gamma (t_n)$ and $\beta_n (1) = \gamma (t_n)$.
The differentiable curves in the family $\mathcal{B}$ fills the gaps along
the trajectory of $\gamma$. Hence, we can  define the parallel transport
along $\gamma$
with respect to $\mathcal{B}$. Precisely, fix a subinterval $[s,t] \subseteq
I$. Take $J = D \cap (s,t]$,
that is, $J$ is the set of the times that jumps occur in $(s,t]$. By
assumption, $J$ is finite and, abusing notation for sack of simplicity,  write
$J = \{t_1 <t_2
< \ldots < t_k\}$.

Now, define the curve $(\gamma \vee \mathcal{B})_{s,t}: [s, t + k] \to M$
which concatenates  $\gamma$ with elements
of the family $\{ \beta_{i}, i=1, \ldots, k \}$ in
the following way:

\[
(\gamma \vee \mathcal{B})_{s,t} (u) =  \begin{cases}
\gamma(u), \ \text{for} \ u \in [s, t_1)  \\
\beta_{1}(u - t_1), \ \text{for} \ u \in [t_1,t_1 + 1] \\
\gamma(u-1), \ \text{for} \ u \in [t_1 + 1, t_2 + 1)  \\
\beta_{2}(u - t_2 - 1), \ \text{for} \ u \in [t_2 + 1, t_2 + 2] \\
\hspace{1cm} \vdots \\
\beta_{k}(u - t_k - k+1), \ \text{for} \ u \in [t_k + k-1, t_k + k] \\
\gamma(u-k), \ \text{for} \ u \in [t_k + k, t + k].  \\
\end{cases}
\]
Since the constructed curve $(\gamma \vee \mathcal{B})_{s,t}$ is continuous
and
differentiable by parts, we define the parallel
transport along $\gamma$ with respect to $\mathcal{B}$  by:
$$P^{\nabla,\mathcal{B}}_{s,t} (\alpha) := P^\nabla_{s,t+k} \Big((\alpha \vee
\mathcal{B})_{s,t}\Big).$$
Note that the domain of $(\gamma \vee \mathcal{B})_{s,t}$ is
artificially extended due to the `fictitious' curves in $\mathcal{B}$ which
fill the gaps of $\gamma$. In next sections, the choice of family
 $\mathcal{B}$ will not be arbitrary. It will be established by the
deterministic
 flow generated by the vector fields of the differential equation.

\bigskip

\section{Stochastic delay differential equations with jumps}

In this section, we present our model of delay differential equations with
jumps (DDEJ), including the deterministic and stochastic case (SDDEJ). The
solution for this
equation is constructed by induction on the number of jumps, such that, after
each jump, we use the theory of differential
equations without jumps. We remark that the existence and
uniqueness of solution in stochastic delay differentiable equations (without
jumps) is a particular case of the theory of stochastic functional differential
equations (see, e.g. Léandre and Mohammed \cite{leandre}).
We start describing the simpler context:

\bigskip

\subsection{Deterministic case}

 Initially,  we construct the discontinuous (c\`adl\`ag) integrator
$S_t$ which drives our model of DDEJ.  Let $(t_n)_{n \in \N}$ be an
increasing, discrete and closed sequence in $\R_{>0}$ which indicates the points
of discontinuities of $S_t$ . Let $(J_n)_{n \in \N}$ be the corresponding
sequence in
$\R$ of the increments at the jumps of $S_t $. Define the integer function
which counts the number of jumps up to time
$t$ by $N_t = \max \{n: t_n \leq t\}$, with the convention that the maximum of
the empty set is zero. Consider the
integrator $S: \R_{\geq 0} \to \R$ in the following way:
$$
S_t = t + \displaystyle \sum_{k=0}^{N_t} J_k.
$$
The DDEJ in the manifold $M$ is written as:
\begin{equation}
\begin{array}{lll} \label{eq: jump-delay}
dx(t) = P^\nabla_{t-d,t} (x) \ F(x(t-d)) \ dS_t
\end{array}
\end{equation}
with initial condition given by a
differentiable curve
$\beta_0: [-d,0] \to M$, where $d \in (0,1]$ is a
fixed time delay and  $F$ is a smooth vector
field in the manifold $M$. We construct a solution $\gamma(t)$ of equation
(\ref{eq:
jump-delay}) as follows:

\bigskip

\noindent $\bullet$ \textbf{Solution before the first jump:}

For $t \in [0, t_1)$, $\gamma(t)$ is the solution of the delay
differential equation (without jumps) given by:

\[
\begin{cases}
x'(t) = P^\nabla_{t-d,t} (x) \ F(x(t-d))  \\
x(t) = \beta_0(t), \ \forall t \in [-d,0].
\end{cases}
\]

\bigskip

\noindent \textbf{$\bullet$ Solution at the jumps:}

Suppose the solution has been constructed in the interval $[0, t_m)$. We define
the solution at the time $t_m$, corresponding to the $m$-th jump. Consider
the ordinary differential equations $y_n'(t) = J_n F(y_n)$, for $n \in
\N$, $n \leq m$. We denote the solution flows of these equations by
$\varphi_t^n$. For each $n\leq m$,  take $z_{n} = \displaystyle
\lim_{s \to t_n-} \gamma(s)$.
Now, let $\mathcal{B}_m = (\beta_n)_{n \leq m}$ be the family of
curves (considered in Section 2) given by:
$$\beta_n (t) = \varphi_t^n (z_{n}) \big|_{[0,1]},$$
and define $\gamma(t_m) = \beta_m(1)$.

\bigskip

\noindent \textbf{$\bullet$ Solution in the intervals between jumps}

In this case, define the solution using the parallel transport along
$\gamma$, with respect to the family $\mathcal{B}_m$. So, for $t \in (t_m,
t_{m+1})$, $\gamma(t)$ is the solution of the
following delay differential equation:

\[
\begin{cases}
x'(t) = P^{\nabla, \mathcal{B}_m}_{t-d,t} (x) \ F(x(t-d))  \\
x(t) = \gamma(t), \ \forall t \in [-d,t_m].
\end{cases}
\]
Note that, although the initial condition of the equation above may have jumps,
results of the standard theory of delay differential equations on existence and
uniqueness still hold. In fact, this initial condition is used only
to parallel transport the vector field, and this (concatenation of) parallel
transport has been defined in the previous section.

\

Therefore, by induction, the solution is well defined for all $t \geq 0$.
Uniqueness comes from the fact that the solution is unique in each step.

The fictitious curves (as in Section 2) we have used in this construction
were established by the deterministic flow of the vector field (without
delay). For future reference, we
call this family of curves associated to $\gamma$ by $\mathcal{B}_F =
\displaystyle \{
\mathcal{B}_n, n \in \N \} $. This model reflects the main motivation that,
in  some physical situations, the informations
arriving at a receptor come from different sources (with corresponding
different delays): here, continuous informations have a fix delay $d$, but
discontinuities in the driver integrator have no delay.

\bigskip

\subsection{Stochastic Case}

Let $(B^1_t, \ldots, B^m_t)$
 be a Brownian motion in a filtered probability space $(\Omega, \mathcal{F},
 \mathcal{F}_t, \Prob )$ and $(N_t)_{t \geq 0}$ be a random counting
process that indicates the number of jumps up to time $t$, with the properties
that $N_0 = 0$ and $N_t$ is finite (almost surely) for all $t \geq 0$.
Consider
 a sequence $(J_k)_{k \in \N}$ of random variables in $\R^{m+1}$. Taking
$B^0_t = t$, the integrator of our model is $L_t= (L^0_t,L^1_t,  \ldots ,
L^{m}_t)$, given by:
$$
L^i_t = \displaystyle B^i_t + \displaystyle \sum_{k=0}^{N_t} J^i_k.
$$
An example of this kind of process is the Levy-jump diffusion (see e.g.
Applebaum \cite{ap}), where $N_t$ is a Poisson process and $(J_k)$ are
i.i.d. random variables.

Write the stochastic delay
differential equation with jumps (SDDEJ) by:
\begin{equation}\begin{array}{lll} \label{eq: jump-delay stochastic}
dx_t = P^\nabla_{t-d,t} A(x_{t-d}) \diamond dL_t,
\end{array}
\end{equation}
where $A_0, A_1, \ldots, A_m$ are smooth vector fields in $M$, considering
initial condition  $\beta_0: [-d,0] \to M$, a differentiable curve in $M$.

We define the solution of this equation in an analogous way to the
deterministic case. So, fixing $\omega \in \Omega$, in the intervals between
the jumps, the solution is given by the corresponding Stratonovich stochastic
delay differential equation, that is:
$$
dx_t = \displaystyle \sum_{i=0}^{m} P^\nabla_{t-d,t} \ A^i (x_{t-d}) \circ
dL_t^i,
$$
with the appropriate initial condition, as we have done in the previous
case. Besides that, at the times of jump, the solution hops instantaneously in
the direction of the solution  at time one of the following ODE
(without delay):

\[
\begin{cases}
y'(t) = \displaystyle \sum^m_{k=0} J^k_n A^k(y) \\
y(0) = \displaystyle \lim_{s \to t_n-} \gamma(s),
\end{cases}
\]

\bigskip

So, with the same notation as before, we have the following result for SDDEJ:

\begin{theorem}
There exists a unique solution $\gamma$ for the SDDEJ (\ref{eq: jump-delay
stochastic}) defined in $t \in [-d, \infty)$, with initial condition
$\gamma(t) = \beta_0(t)$ when $t\in [-d,0]$.

\begin{proof}
The existence follows by an analogous construction we have done for the
deterministic case. Uniqueness holds since that, in each
step of the construction, the respective solution is unique, by theory of ODE
and standard stochastic delay differential equations, as in \cite{leandre}.
\end{proof}

\end{theorem}

In this case, for each $\omega \in \Omega$, we have a family of
differentiable curves established by the deterministic flow of the vector 
field. Again, we call this random family of jumps associated to $\gamma$ by 
$\mathcal{B}_F (\omega) = \displaystyle \{
\mathcal{B}_n (\omega), n \in \N \} $.

\medskip

The idea of jumping in the direction of the deterministic flow at time one
comes from Marcus SDE, where the integrator is a semimartingale with jumps (for more details see
Kurtz, Pardoux e Protter \cite{KPP}).

\bigskip

\section{Geometrical aspects of SDDEJ}

In this section, we show that the parallel transport, i.e., the horizontal lift
of a solution of an SDDEJ in a manifold $(M, \nabla)$ can be described as an
SDDEJ in the linear frame bundle $BM$, with respect to a horizontal connection
in $BM$ (described below). The equation for this
horizontal lift corresponds to an extension of the results on stochastic
geometry started with It\^o \cite{Ito} and Dynkin \cite{Dynkin} to our
model of SDDEJ (see also \cite{PP}).

\subsection{Horizontal lifts to the frame bundle}

For reader's convenience we recall briefly some geometrical facts about the
frame bundle of a
manifold (for more details, see e.g., among many others, Elworthy
\cite{elworthy},  Kobayashi
and Nomizu \cite{nomizu}). Let $M$ be a differentiable manifold, with dimension
$n$. The frame bundle $BM$ of $M$ is the set of all linear
isomorphisms $p: R^n \to T_x M$ for $x \in M$. The projection $\pi: BM \to
M$ maps $p$ to the corresponding  $x\in M$. $BM$ is a principal bundle over
$M$, with
right action of the Lie group $GL(n,\R)$, given by the composition with the
linear isomorphisms.

\medskip

Given $p$ in the manifold $BM$, each tangent space $T_p BM$ can be decomposed
as a direct sum
of the vertical and a horizontal subspace, $T_p BM=V_p BM \oplus H_p BM $. The
vertical subspace is determined by
$V_p BM = Ker (\pi_* (p))$, where $\pi_*$ denotes the derivative of the
projection $\pi$. The horizontal subspace $H_p BM$ is established
by the connection $\nabla$ in $M$, namely it is generated by the derivative of
parallel frames along curves in $M$ passing at $\pi (p)$. In this context,
one can consider the horizontal lift of a vector $v \in T_x M$ at $p \in
\pi^{-1}(x)$ as the unique tangent vector $v^H \in H_p BM$ such that $\pi_* (p)
(v^H) = v$.

\medskip

We say that a differentiable curve $\alpha: I \to BM$ is horizontal when its
derivative belongs to $H_{\alpha(t)} BM$ for all $t \in I$. In fact, given a
differentiable curve $\beta: [0,T) \to M$, and $p \in \pi^{-1}
(\beta(0))$, there exists a unique horizontal curve $\beta^H: [0,T) \to BM$,
with the property that $\pi(\beta^H(t)) = \beta(t)$ for all $t$ in the domain.
The curve $\beta^H$ is called the horizontal lift of the curve $\beta$
(see e.g. \cite{nomizu}).

\medskip

Putting together the technique of Section 2 and the horizontal lift described
above, one can define the horizontal lift of a curve with jumps in $M$. Let
$\gamma: [0,\infty) \to M$ be a càdlàg curve, $D = \{t_1, t_2,
\ldots \}$ the countable, closed and discrete set of points of discontinuity, and
$\mathcal{B} = (\beta_n)_{n \in \N}$ be a family of differentiable curves
$\beta_n: [0,1] \to M$, such that, for all $n \in \N$, $\beta_n (0) =
\displaystyle \lim_{s \to t_n-} \gamma (t_n)$ and $\beta_n (1) = \gamma (t_n)$.
Fix $p \in \pi^{-1}(0)$. Under these conditions, we define the horizontal lift
of $\gamma$ in $p$ with respect to the family $\mathcal{B}$ by the càdlàg curve
in $BM$:
$$
\gamma^{H,\mathcal{B}}_p (t) := P^{\nabla,\mathcal{B}}_{0,t} (\gamma) \circ
p.
$$

Each element $A$ in the Lie algebra $\mathcal{G}l(n,\mathbb{R})$
of the Lie
group $GL(n,\R)$ determines a vertical vector field in $BM$ given by, at
a point $p\in BM$,
\[
A^{*}(p) = \frac{d}{dt}\ \big( p \cdot  \exp{At} \big) |_{t=0}.
\]
The map $\mathcal{G}l(n,\mathbb{R}) \mapsto V_pBM$ is surjective.
In order to define a SDDEJ in $BM$, one needs a connection in this
manifold as well. There are many ways of extending a connection $\nabla$ of $M$
to $BM$. In this section we are interested in the so called horizontal lift
$\nabla^{H}$ which is defined
(for a torsion free connection $\nabla$, see e.g. Cordero et al. \cite[Chap.
6]{cordero}) as the unique connection in $BM$
which satisfies:
\begin{equation}  \label{conexao-horizontal}
 \left\{ \begin{array}{lcl}
 \vspace{1.72mm} \nabla^{H}_{A^{*}}B^{*} & = & (AB)^{*} \\
 \vspace{1.72mm} \nabla^{H}_{A^{*}}X^{H} & = & 0  \\
 \vspace{1.72mm}  \nabla^{H}_{X^{H}}A^{*} & = & 0  \\
 \nabla^{H}_{X^{H}}Y^{H} & = & (\nabla_{X}Y)^{H}
\end{array} \right.
\end{equation}
 This extension has the property that parallel transport
commutes with the horizontal lift, that is, if $\alpha$ is a curve in $BM$, for
any $v \in T_{\pi \circ \alpha(0)}M$, it holds that
$P^{\nabla^{H}}_{0,t}(\alpha)(v^{H}) = (P^{\nabla}_{0,t}(\pi \circ
\alpha)(v))^{H}$, see \cite[Lemma 2.1]{PP}.

\subsection{Main results}

As in the previous section, initially we deal with deterministic systems.
Catuogno and
Ruffino \cite[Prop. 2.1]{PP} used the extended connection $\nabla^H$ to
prove that, if a curve $\gamma$ is the solution of a deterministic DDE in $M$,
then its horizontal lift is the solution of the corresponding DDE in $BM$.
Now, we extend this result to DDEJ. Let
$F: M \to TM$ be a smooth vector field and $\beta_0: [-d,0] \to M$ a
differentiable curve (initial condition).  Let $\gamma$ be the solution of
the following delay differential
equation with jumps in $M$:

\[ \begin{cases} dx(t) = P^{\nabla}_{t-d,t} (x) \ F(x(t-d)) \ dS_t \\
x(t) = \beta_0 (t) \ \ \text{for} \ t \in [-d,0]. \end{cases}
\]
This solution induces canonically the family $\mathcal{B}_F$ of
curves along the deterministic flow at the jumps, as defined in Section
3.1. Fix $p \in \pi^{-1} (\gamma(-d))$. Consider $\beta_0^H: [-d,0] \to BM$ the horizontal lift of $\beta_0$ at the point $p$, and let $F^H$ be the horizontal lift of the vector field $F$.
Under these conditions:

\bigskip

\begin{theorem}
\label{lift}
The horizontal lift $\gamma^{H,\mathcal{B}_F}_p$ is the solution of the following DDEJ in $BM$ with respect to the connection $\nabla^H$:

\begin{equation}\begin{array}{lll} \label{eq: lifting}
\begin{cases} du(t) = P^{\nabla^H}_{t-d,t} (u) \ F^H(u(t-d)) \ dS_t \\
u(t) = \beta^H_0 (t) \ \ \text{for} \ t \in [-d,0]. \end{cases}
\end{array}
\end{equation}

\begin{proof} Let $(t_n)_{n \in \N}$ be the increasing sequence of
discontinuities of the integrator $S_t$. Let $u: [0, \infty) \to BM$ be the 
solution of equation (\ref{eq: lifting}), whose existence and uniqueness is 
guaranteed in Section 3.1. We show, by induction on the number of jumps, that 
$u(t) = \gamma^{H,\mathcal{B}_F}_p (t)$ for all $t \geq 0$.

For $t \in [0,t_1)$, the solution $u(t)$ is given by a delay
differential equation, so we can apply the result
without jumps as in  \cite[Prop. 2.1]{PP} to obtain the equality. In 
particular, define $p_1:=\displaystyle \lim_{s \to {t_1}^-} u(s) = 
\displaystyle \lim_{s \to
{t_1}^-} \gamma^{H,\mathcal{B}_F}_p (s).$

\medskip

At time $t_1$, when the first jump occurs, take the curve $\beta_1 \in
\mathcal{B}_F$ (we recall that this is the solution of 
the
ODE $y'(t) = J_1 F(y)$, with initial condition $y(0)
= \pi(p_1)$). Consider its horizontal lift $\beta^H_1$  at $p_1$. As $\beta^H_1$ 
is the solution of $z'(t) = J_1 F^H (z)$, with initial condition $z(0) = p_1$,
we have that $u(t_1) = \gamma^{H,\mathcal{B}_F}_p (t_1)$.
Now, arguing by induction, suppose that $u(t) = 
\gamma^{H,\mathcal{B}_F}_p (t)$ for all $t \in [-d, t_m]$ for  $m \geq 1$. We 
claim that this equality also holds in the interval $(t_m, t_{m+1}]$. In fact, 
let $k$ be the number of jumps that occur in the interval $(t_m - d, t_m)$.

\medskip

Firstly, for the case $(t_{m+1} - t_m) > d$, we have:
$$t_m < t_{m-k} + d < \ldots < t_{m-1} + d < t_m + d < t_{m+1},$$
and it is enough to analyse each of these $(k+ 2)$ subintervals. For the first 
subinterval, that
is, $(t_m, t_{m-k} + d)$, we consider a DDE with delay $(d+k+1)$, where the
initial condition is the concatenation of the following curves:
\begin{itemize}
\item $u_t$, in the interval $[t_m - d, t_{m-k})$;
\item $\beta^H_{m-k}$, in the interval [0,1];
\item $u_t$, in the interval $[t_{m-k}, t_{m-k+1})$;
\item $\beta^H_{m-k+1}$, in the interval [0,1];
\

\hspace{1cm} \vdots
\item $u_t$, in the interval $[t_{m-1}, t_m)$;
\item $\beta^H_{m}$, in the interval [0,1].
\end{itemize}

Therefore, for $t \in (t_m, t_{m-k} + d)$, we have a delay differential
equation, with a fictitious bigger delay, but without jumps. Applying the 
result in \cite[Prop. 2.1]{PP}
the equality in this subinterval holds. For the $n$-th subinterval, with 
$n<(k+2)$, consider a DDE (without jumps), where the initial condition 
concatenates appropriately $u_t$ and $\beta^H_{i}$, with $i \in \{m-k+n-1, 
\ldots, m\}$ and use \cite[Prop. 2.1]{PP}. And for the last subinterval, that
is, for $t\in (t_m + d, t_{m+1})$, take the DDE (without jumps) with delay $d$ 
and
initial condition given by $u_t$, and again we have the equality.

\medskip

Secondly, for the case $(t_{m+1} - t_m) \leq d$ the argument 
is essentially the same, just we have to consider a smaller number of 
subintervals. Therefore, we have that $u(t) = 
\gamma^{H,\mathcal{B}_F}_p (t)$ for $t \in (t_m, t_{m+1})$. In particular, 
define:
$$
p_{m+1}:= \displaystyle \lim_{s \to t_{m+1}^-} u(s)  = \displaystyle 
\lim_{s \to t_{m+1}^-} \gamma^{H,\mathcal{B}_F}_p(s).
$$


Finally, consider $\beta^H_{m+1}$, the horizontal lift of the differentiable 
curve $\beta_{m+1}$  at the point $p_{m+1}$. The equality at $t_{m+1}$ also 
holds: $u(t_{m+1}) = 
\gamma^{H,\mathcal{B}_F}_p (t_{m+1})$. The proof is complete.

\end{proof}
\end{theorem}

\medskip
The same result holds in the stochastic case thanks to the transfer 
principle (see, e.g. Emery \cite{emery}), in the sense that: with the same 
notation
as before, let $\alpha$ be the solution of the SDDEJ

\[ \begin{cases} dx(t) = P^{\nabla}_{t-d,t} (x) \ F(x(t-d)) \diamond dL_t \\
x(t) = \beta_0 (t) \ \ \text{for} \ t \in [-d,0]. \end{cases}
\]
This solution induces  the random family $\mathcal{B}_F (\omega)$ of
curves along the deterministic flow at the jumps, as defined in Section 3.2. 
Hence, we have the following:

\begin{theorem}

The horizontal lift $\alpha^{H,\mathcal{B}_F (\omega)}_p$ of $\alpha$ (solution of the SDDEJ in $M$) is the solution of the following SDDEJ in $BM$, with
respect to the connection $\nabla^H$:

\begin{equation}\begin{array}{lll} \label{eq: stochastic lifting}
\begin{cases} du(t) = P^{\nabla^H}_{t-d,t} (u) \ F^H(u(t-d)) \diamond dL_t \\
u(t) = \beta^H_0 (t) \ \ \text{for} \ t \in [-d,0]. \end{cases}
\end{array}
\end{equation}

\begin{proof}
In the proof of the deterministic case (Theorem 4.1), in each step, the
solution is given by a standard delay differential equation (without jumps).
 The corresponding lift is the solution of the lifted equation in $BM$.  So, 
applying successively the transfer
principle: before the first jump, at the jumps and in the intervals between
jumps, we have the result (cf. \cite[Thm. 2.2]{PP}).

\end{proof}
\end{theorem}

\noindent The results above do not exhaust the subject; in fact, it opens the 
possibility of exploring the geometry for SDDEJ and 
corresponding applications: say, holonomies, invariant measures, stability 
(Lyapunov exponents), rotation numbers, etc.


\end{document}